\documentclass[preprint,12pt]{elsarticle}

\usepackage{graphicx}
\usepackage{amssymb}
\usepackage{amsthm}
\usepackage[fleqn]{amsmath}
\usepackage{mathtools}
\usepackage{cases}
\usepackage{lipsum}
\usepackage{leftidx}
\usepackage{float}
\makeatletter
\def\ps@pprintTitle{%
 \let\@oddhead\@empty
 \let\@evenhead\@empty
 \def\@oddfoot{}%
 \let\@evenfoot\@oddfoot}
\makeatother

\biboptions{sort&compress}
\newtheorem{thm}{Theorem}

\theoremstyle{example}
\newtheorem{exm}{Example}

\theoremstyle{definition}
\newtheorem{defn}{Definition}
\theoremstyle{remark}
\newtheorem{rem}{Remark}

\journal{}

\begin{document}

\begin{frontmatter}



\title{The flaw in the conformable calculus:
it is conformable because it is not fractional}
\author{Ahmed A. Abdelhakim}
\address{Mathematics Department, Faculty of Science, Assiut University, Assiut 71516 - Egypt\\
Email: ahmed.abdelhakim@aun.edu.eg}
\begin{abstract}
We point out a major flaw in the conformable calculus. We demonstrate why it fails at defining a fractional derivative and where exactly these tempting conformability properties come from.
\end{abstract}
\end{frontmatter}
\section{Introduction}
\indent Khalil et al. proposed a definition for the fractional derivative in \cite{khalil} using what they called the ``conformable derivative". This concept was quickly adopted by T. Abdeljawad in \cite{Abdeljawad} where he
claims to have developed some tools of
fractional calculus.
\\
\indent The conformable derivative is local
by its very definition. Moreover, we proved rigorously in \cite{ahmedMachado} that the conformable derivative of a function $f$ does not exist at any point $x>0$, unless $f$ is differentiable at $x$.
The term ``conformable"
is supposedly attributed to the properties
this proposed definition provides. \\
\indent We point out the flaw in Khalil et al.'s definition and uncover
the real source of this conformability through reviewing the statements and proofs in \cite{Abdeljawad,khalil}.
Analogous remarks apply to the statements and proofs in \cite{Anderson} and
\cite{Katugampola1}. It turns out that the reason behind the conformability of this derivative
is, ironically, the same reason it is not fractional.\\
\indent We would like to emphasize here that we are not reviewing the aforementioned work to provide useful formulae to work with. On the contrary,
our real purpose is to discourage researchers from using it, by making it clear from the mathematical point of view why the conformable derivative is not fractional.\\
\indent We have shown in (\cite{ahmedMachado}, Section 5) the disadvantages of using the conformable definition in solving fractional differential equations. It breaks the fractional equation and replaces it with an ordinary equation that may no longer properly describe the underlying fractional phenomenon. This is probably the reason it produces a substantially larger error compared with the Caputo fractional derivative when used to solve fractional models (see \cite{ahmedMachado}, Section 6).\\
\indent
We discuss concrete examples
that illustrate how the conformable derivative
is incapable of giving the fractional derivative
obtainable from the classical Riemann-Liouville
or Caputo derivatives. More examples are provided
to show how the conformable operator produces
functions with a much different behaviour than
the classical fractional derivatives. The latter
are known to be successful at describing many fractional phenomena (see e.g. \cite{Mainardi}).
\section{The problems in the statements and proofs in \cite{khalil}}
The results in \cite{khalil}
are all based on the following definition:
\begin{defn}\label{defn1}(\cite{khalil}, Definition 2.1)
Given a function
$f : [0,\infty[ \rightarrow \mathbb{R}$, then the `` conformable fractional derivative" of $f$ of order
$\alpha$ is defined by
\begin{equation}\label{q011}
T_{\alpha}f(t):=
\lim_{\epsilon\rightarrow 0}
\frac{f(t+\epsilon t^{1-\alpha})-f(t)}{\epsilon},
\end{equation}
for all $t>0$, $\alpha\in\, ]0,1[$.
\end{defn}
Definition \ref{defn1} is flawed. Once we establish this, it will be immediately seen that the proofs in \cite{khalil} are unnecessarily involved. More importantly, the results therein will be found insignificant
as they follow directly from the traditional integer-order calculus.\\
\indent
We proved the following theorem in \cite{ahmedMachado}:
\begin{thm}\label{thm1}(\cite{ahmedMachado}, Theorem 1) Fix $\,0<\alpha<1$ and let $t>0$. A function
$\,f:[0,\infty[\,\longrightarrow \mathbb{R}\,$
has a ``conformable fractional derivative" of
order $\alpha$ at $t$ if and only if
it is differentiable at $t$, in which case we have the pointwise relation
\begin{equation}\label{pntwserltn1}
T_{\alpha}f(t)=
t^{1-\alpha}f^{\prime}(t).\end{equation}
\end{thm}
We note here the problems
with Definition \ref{defn1}
in the light of Theorem \ref{thm1}:
\begin{rem}
 The limit (\ref{q011}) does not exist unless $\displaystyle\lim_{\epsilon\rightarrow 0}
  {\left(f(t+\epsilon)-f(t)\right)}/{\epsilon}$ exists.
In other words, there does not exist a function
  differentiable in the sense of
 Definition \ref{defn1} that is not differentiable.
In fact, the false claim in
\cite{Abdeljawad,Anderson,khalil,Katugampola1}
that the  ``conformable" derivative may exist at a point where the function is not differentiable is the only excuse for the results in these papers.
\end{rem}
\begin{rem}\label{smconf}
The identity (\ref{pntwserltn1}) is the reason $T_{\alpha}f$ demonstrates ``conformability". The conformability comes precisely from the integer-order derivative,
      the factor $f^{\prime}$, in (\ref{pntwserltn1}).
\end{rem}
\begin{rem}\label{smfrac}
The derivative $T_{\alpha}f$ is not a fractional derivative. It is exactly the integer-order derivative times the root function $t^{1-\alpha}$.
\end{rem}
Therefore, Definition \ref{defn1} is to be understood as follows:
\begin{defn}\label{correctdef}(What Definition 2.1 in \cite{khalil} really suggests)
Given a function
$f : [0,\infty[ \rightarrow \mathbb{R}$, then
$f$ is $\alpha$-differentiable at $t>0$,
if it is differentiable at $t$, and
its $\alpha$-derivative
$\,T_{\alpha}f(t):=
t^{1-\alpha}f^{\prime}(t),\; t>0$.
\end{defn}
Now, we show how this correct understanding of what Definition \ref{defn1} proposes trivializes the results in \cite{khalil}.
\begin{thm}\label{khthm1}(\cite{khalil}, Theorem 2.1)
If a function
$f : [0,\infty[ \rightarrow \mathbb{R}$ is $\alpha$-differentiable at $t_{0}>0$, $\alpha\in\, ]0,1]$, then $f$ is continuous at $t_{0}$.
\end{thm}
If $f$ is $\alpha$-differentiable at $t_{0}>0$, then
it is differentiable at $t_{0}$.
It is well-known that if a function is differentiable at some point, then it is continuous thereat.\\
The next theorem explains why $T_{\alpha}f$ is described as conformable. We show why the statements are trivial
and how the conformability comes from (\ref{pntwserltn1}).
\begin{thm}\label{khthm2}(\cite{khalil}, Theorem 2.2)
Let $\alpha\in\,]0,1]$ and $f,g$ be
$\alpha$-differentiable. Then
\begin{description}
 \item[\textbf{(1)}] $T_{\alpha}(af+bg)=aT_{\alpha}(f)+bT_{\alpha}(g)$
     for all $\;a,b\in \mathbb{R}$.
\item[\textbf{(2)}]$T_{\alpha}(t^{p})=pt^{p-\alpha}$
for all $p\in \mathbb{R}$.
\item[\textbf{(3)}] $T_{\alpha}(f)=0$ for all constant
functions $f$.
\item[\textbf{(4)}]$T_{\alpha}(fg)=
fT_{\alpha}(g)+gT_{\alpha}(f)$.
\item[\textbf{(5)}]$T_{\alpha}\left(\frac{f}{g}\right)=
\frac{gT_{\alpha}(f)-fT_{\alpha}(g)}{g^2}$.
\item[\textbf{(6)}]
If, in addition, $f$ is differentiable, then
$T_{\alpha}(f)(t)=t^{1-\alpha}f^{\prime}(t)$.
\end{description}
\end{thm}
If $f,g$ are $\alpha$-differentiable, then they are in fact differentiable, and we have
$T_{\alpha}(f)(t)=t^{1-\alpha}f^{\prime}(t)$, and
$T_{\alpha}(g)(t)=t^{1-\alpha}g^{\prime}(t)$.\\
\indent Let us start with \textbf{(1)}. Since $f,g$ are differentiable, then so is $af+bg$. By Theorem \ref{thm1},
$af+bg$ is $\alpha$-differentiable and we have
\begin{equation*}
T_{\alpha}(af+bg)=
t^{1-\alpha}(af+bg)^{\prime}=
a t^{1-\alpha}f^{\prime}+
b t^{1-\alpha}g^{\prime}=aT_{\alpha}(f)+bT_{\alpha}(g).
\end{equation*}
The proofs of items \textbf{(2)} through \textbf{(5)} are as trivial as the proof of \textbf{(1)}.\\
\indent The statement \textbf{(6)} is inaccurate.
The truth is $f$ is $\alpha$-differentiable at $t>0$ if and only if $f$ is differentiable at $t$. Thus, if
f is $\alpha$-differentiable at $t>0$, then
$T_{\alpha}(f)(t)=t^{1-\alpha}f^{\prime}(t)$.
We do not need to require $f$ to be differentiable. Differentiability is already implied by assuming $f$ is $\alpha$-differentiable.
\begin{thm}\label{khthm3}(\cite{khalil}, Theorem 2.3)
Let $a>0$ and $f:[a,b]\rightarrow \mathbb{R}$ be a given function such that:
\begin{description}
  \item[(i)] $f$ is continuous on $[a, b]$,
  \item[(ii)] $f$ is $\alpha$-differentiable for some $\alpha \in\,]0,1[$,
  \item[(iii)] $f(a)=f(b)$.
\end{description}
Then, there exists $c \in\,]0,1[$ such that
$ f^{(\alpha)}(c) = 0$.
\end{thm}
The condition (ii) implies that $f$
is differentiable on $]a,b[$
and $f^{(\alpha)}(t)=t^{1-\alpha}
f^{\prime}(t)$ for all $t\in\,]a,b[$. We know from the classical Rolle's theorem that there exists $c\in\,]a,b[$ such that $f^{(\alpha)}(c)=c^{1-\alpha}f^{\prime}(c)=0$.
\begin{thm}\label{khthm4}(\cite{khalil}, Theorem 2.4)
Let $a>0$ and $f:[a,b]\rightarrow \mathbb{R}$ satisfy
\begin{description}
  \item[(i)] $f$ is continuous on $[a, b]$,
  \item[(ii)] $f$ is $\alpha$-differentiable for some $\alpha \in\,]0,1[$.
\end{description}
Then, there exists $c \in\,]0,1[$ such that
$ f^{(\alpha)}(c) = \frac{f(b)-f(a)}
{\frac{1}{\alpha}b^{\alpha}-\frac{1}{\alpha}a^{\alpha}}$.
\end{thm}
Once again, by Theorem \ref{thm1}, the condition
(ii) implies that $f$ is differentiable
on $]a,b[$ and $f^{(\alpha)}=t^{1-\alpha}
f^{\prime}$ on $]a,b[$. Now, apply
the classical Cauchy mean value theorem to the functions $f$ and $t\mapsto \frac{t^{\alpha}}{\alpha}$ on the interval $]a,b[$, we already know that there is $c\in
\,]a,b[$ such that
\begin{equation*}
f^{(\alpha)}(c)=
\frac{f^{\prime}(c)}{c^{\alpha-1}}=
\frac{f(b)-f(a)}
{\frac{1}{\alpha}b^{\alpha}-\frac{1}{\alpha}a^{\alpha}}.
\end{equation*}
\indent Let $a>0$. Proposition 2.1 in \cite{khalil} introduces absolutely no novelty because, by (\ref{pntwserltn1}),
and the fact that $t\mapsto t^{1-\alpha}$
is locally bounded, $f^{(\alpha)}$ is bounded
on $[a,b]$ if and only if
$f^{\prime}$ is bounded on $[a,b]$. And if
$f^{\prime}$ is bounded on $[a,b]$, then
$f$ is Lipschitz on $[a,b]$, not only uniformly continuous. \\
\indent The remark that follows Proposition 2.1 in \cite{khalil} is false. If an $\alpha$-differentiable function $f$ on $]a,b[$ is uniformly continuous on $[a,b]$, then its $\alpha$-derivative is not necessarily bounded therein. A counterexample is $f(t)=t^{\frac{1}{4}}$
which is Lipschitz on $[0,1]$, but $f^{\frac{1}{2}}(t)=
t^{-\frac{1}{4}}/4$ is unbounded on $]0,1[$.\\
\indent Take a look at
\begin{defn}\label{defn2}(\cite{khalil}, Definition 3.1) Let $a\geq 0$. The $\alpha$-integral
of a function $f$ is
$I^{a}_{\alpha}(f)(t):=
\int_{a}^{t}\frac{f(x)}{x^{1-\alpha}}dx$.
\end{defn}
The example given in \cite{khalil},
right after Definition 3.1, seems to try to sell
$I_{\alpha}^{a}$ as the antiderivative of
$T_{\alpha}$. Of course, $T_{\frac{1}{2}}\left(
\sin{t}\right)=
\sqrt{t}\cos{t}$, and
$I^{0}_{\frac{1}{2}}\left(\sqrt{t}\cos{t}\right)=
\int_{0}^{t}\cos{x}dx=\sin{t}$.
The truth is
$\int_{a}^{t}\frac{f^{(\alpha)}(x)}{x^{1-\alpha}}dx=
\int_{a}^{t} f^{\prime}(x)dx=f(t)-f(a)$.
For example,
$I^{0}_{\frac{1}{2}}\left(\sqrt{t}\sin{t}\right)=
I^{0}_{\frac{1}{2}}\left(T_{\frac{1}{2}}(-\cos{t})\right)=
\int_{0}^{t}\sin{x}dx=1-\cos{t}$.
\section{The problems in the statements and proofs in \cite{Abdeljawad}}
We proceed to demonstrate the
flaws in the definitions suggested in  \cite{Abdeljawad}. We prove that
the tools of calculus proposed there lack the novelty, as they are trivial consequences of
the traditional calculus. The ideas in \cite{Abdeljawad} are all based on the following definition:
\begin{defn}\label{defn2}(\cite{Abdeljawad}, Definition 2.1)
The (left) fractional derivative starting from $a$ of a function $f:[a,\infty[\rightarrow \infty$ of order
$0<\alpha\leq 1$ is defined by
\begin{equation}\label{q022}
T^{a}_{\alpha}f(t):=
\lim_{\epsilon\rightarrow 0}
\frac{f(t+\epsilon (t-a)^{1-\alpha})-f(t)}{\epsilon},\;t>a.
\end{equation}
The (right) fractional derivative of order $0<\alpha\leq 1$ of a function $f:]-\infty,b ]\rightarrow \infty$ is defined by
\begin{equation}\label{q023}
^{b}_{\alpha} Tf(t):=-
\lim_{\epsilon\rightarrow 0}
\frac{f(t+\epsilon (b-t)^{1-\alpha})-f(t)}{\epsilon},\;
t<b.
\end{equation}
If $T^{a}_{\alpha}f(t)$ exists on
$]a,b[$ then $T^{a}_{\alpha}f(a):=
\lim_{t\rightarrow a^{+}}T_{\alpha}f(t)$.
If $^{b}_{\alpha}Tf(t)$ exists on
$]a,b[$ then
$\, {^{b}_{\alpha} }Tf(b):=
\lim_{t\rightarrow b^{-}} {^{b}_{\alpha}Tf(t)}$.
\end{defn}
\indent It is also noted in \cite{Abdeljawad}
that if $f$ is differentiable, then
$\,T^{a}_{\alpha}f(t)= (t - a)^{1-\alpha} f^{\prime}(t)$ and $^{b}_{\alpha} Tf(t)= -(b - t)^{1-\alpha}f^{\prime}(t)$. \\\\
\indent The following theorem is given in \cite{ahmedMachado}:
\begin{thm}\label{thm02}(\cite{ahmedMachado}, Theorem 3) Suppose $h:\,]-1,1[\,\times\mathbb{R}\longrightarrow
\mathbb{R}$ is such that
$\,\lim_{\epsilon\rightarrow 0}h(\epsilon,t_{0})\neq 0$
for some $t_{0}\in \mathbb{R}$. Then a function $\psi:\mathbb{R}\longrightarrow
\mathbb{R}$ is differentiable at $t_{0}$
if and only if the limit
\begin{equation*}
\tilde{\psi}(t_{0}):=\lim_{\epsilon\rightarrow 0}
\frac{\psi\left(t_{0}+\epsilon h(\epsilon,t_{0})\right)-\psi(t_{0})}{\epsilon}
\end{equation*}
exists, in which case $\tilde{\psi}(t_{0})=
\psi(t_{0})\psi^{\prime}(t_{0})$,
$\; \psi(t)=
\lim_{\epsilon\rightarrow 0}h(\epsilon,t)$.
\end{thm}
\indent Let us see the problems
in Definition \ref{defn2}:
\begin{rem}\label{remdfn}
According to Theorem \ref{thm02}, the limit in (\ref{q022}) exists at $t>a$
if and only if  $\displaystyle\lim_{\epsilon\rightarrow 0}  {\left(f(t+\epsilon)-f(t)\right)}/{\epsilon}$. exists. Similarly, the limit in (\ref{q023}) exists at $t<b$ if and only if $f^{\prime}(t)$ exits.
This means that neither $T^{a}_{\alpha}f(t)$
nor $^{b}_{\alpha} Tf(t)$ exits unless
$f$  is differentiable at $t$. In fact, by Theorem
\ref{thm02}, Definition \ref{defn2} reads:
\begin{eqnarray}\label{0qq3}
\begin{split}
  T^{a}_{\alpha}f(t)&:=& (t-a)^{1-\alpha}
  f^{\prime}(t),\; t>a \\
^{b}_{\alpha} Tf(t)&:=& (b-t)^{1-\alpha}
  f^{\prime}(t),\;t<b,
  \end{split}
\end{eqnarray}
provided $f^{\prime}(t)$ exists.
\end{rem}
\begin{rem}
Unlike with the classical fractional derivatives of Riemann-Liouville and Caputo,
Definition \ref{defn2} does not work for functions defined on $\mathbb{R}$. Indeed, by Remark \ref{remdfn}, if $f$ is defined on $\mathbb{R}$, then both derivatives
$T^{-\infty}_{\alpha}f(t)$ and
$^{\infty}_{\alpha}Tf(t)$ are ill-defined
at every $t\in \mathbb{R}$, which is unacceptable.
\end{rem}
\begin{rem}
There is no geometric or physical motivation
that justifies the negative sign
in the definition of the operator
$^b_{\alpha}T$.
Furthermore,
the case $\alpha=1$ is supposed to give the
left first order
integer derivative, but Remark \ref{remdfn} implies
\begin{equation*}
^b_{1}Tf(t)=-f^{\prime}(t),\;t<b,
\end{equation*}
which is neither the left nor the right derivative
of $f$ at $t$.
\end{rem}
\begin{rem}
There is an obvious inconsistency in
Definition \ref{defn2} when it comes to defining
$T^{a}_{\alpha}f(a)$ and $^{b}_{\alpha}Tf(b)$.
Let $a<b$. We have clarified in Remark \ref{remdfn} that if a function is not differentiable
at some point $t\in\,]a,b[$, then
the pointwise criterion of Definition \ref{defn2} does not allow it to be $\alpha$-differentiable at $t$.
This, however, excludes the endpoints $a$ and $b$.
We are going to discuss
$T^{a}_{\alpha}f(a)$ and the analogue applies to
$^{b}_{\alpha}Tf(b)$.
Unjustifiably, Definition \ref{defn2} allows
the derivative $T^{a}_{\alpha}f(a)$
to exist regardless of the existence
of the right derivative of $f$ at $a$.
Precisely, by Remark \ref{remdfn},
if $T^{a}_{\alpha}f$ exists on
$]a,a+\delta[$, for some $\delta>0$, then
\begin{equation*}
T^{a}_{\alpha}f(a)=
\lim_{t\rightarrow a^{+}}
T^{a}_{\alpha}f(t)=
\lim_{t\rightarrow a^{+}}
(t-a)^{1-\alpha}f^{\prime}(t).
\end{equation*}
Therefore, according to Definition \ref{defn2},
 $T^{a}_{\alpha}f(a)$ exists
if and only if
$f^{\prime}$ exists on $]a,a+\delta[$, and $\lim_{t\rightarrow a^{+}}
(t-a)^{1-\alpha}f^{\prime}(t)$ exists.
This is evidently a weaker condition than the existence
of $f^{\prime}$ on $]a,a+\delta[$
and $\lim_{t\rightarrow a^{+}}f^{\prime}(t)$.
It is also independent of the existence
of the right derivative $f^{\prime}_{+}(a)$ of  $f$
at $a$. Many examples are given in \cite{Abdeljawad}
for functions differentiable on $]a,b[$ such that $T^{a}_{\alpha}f(a)$
exists, but $f^{\prime}_{+}(a)$ does not. This may lead to the false intuition that
the operator $T^{a}_{\alpha}$ is well-defined on a larger class of functions than the derivative.  The reality is there exist smooth functions
on $]a,b[$ such that $f^{\prime}(a)$ exists but
$T^{a}_{\alpha}f(a)$ does not. Consider for instance
$\displaystyle g(t):=
\left\{
    \begin{array}{ll}
       x^2\sin{\frac{1}{x^3}}, & \hbox{$x\neq 0$;} \\
     0, & \hbox{x=0.}
     \end{array}
      \right.$. We have $g\in C^{\infty}(\mathbb{R}\setminus\{0\})$ and $g^{\prime}(0)=0$, yet
$\lim_{x\rightarrow 0^{+}}
x^{1-\alpha}g^{\prime}(x)$ does not exist for any
$0\leq\alpha\leq 1$.
\\
\indent
Another issue with $T^{a}_{\alpha}f(a)$
is that its existence depends
on the domain. For example,
$h(x):=\sin{\sqrt{t-t_{0}}}$ is not differentiable at $t=t_{0}$, and consequently, by Remark \ref{remdfn},
$T^{c}_{\alpha}h(t_{0})$ does not exist. But this is true only if $h$ is considered on the domain $[c,\infty[$ with any $c<t_{0}$. If $c=t_{0}$, however, then  $T^{c}_{\alpha}h(t_{0})$ magically exists and equals 0,
for every $0<\alpha\leq \frac{1}{2}$.
\end{rem}
\begin{rem}\label{remdfn1}
The identities (\ref{0qq3}) prove that
the derivative in Definition \ref{defn2}
is not fractional and that the conformability comes
from the integer-order derivative factor.
What is worse is that the derivative in Definition \ref{defn2} fails
to give the fractional derivative for some functions
whose fractional derivative exist and can be
easily calculated using the Riemann-Liouville or Caputo definition. See the following examples:
\begin{exm}
Consider the function
$f_{1}(t):=\chi_{\raisebox{-.5ex}{$\scriptstyle [0,1]$}}{(t)}$ on $[0,2]$. It is easily verifiable that $(T^{0}_{\alpha}f_{1})(1)$ does not exist. But the Riemann-Liouville fractional derivative
$D^{\alpha}_{0^{+}}f_{1}$
exists at $t=1$, and
$\, D^{\alpha}_{0^{+}}f_{1}(1)=\frac{1}{\Gamma(1-\alpha)}
\int_{0}^{1}\frac{d\xi}{(1-\xi)^{\alpha}}=
\frac{1}{(1-\alpha)\Gamma(1-\alpha)}$.
\end{exm}
\begin{exm}
Consider the function
$f_{2}(t):=|t-1|$ on $[0,2]$. Again, $(T^{0}_{\alpha}f_{2})(1)$ does not exist. Nevertheless, the Caputo fractional derivative $\leftidx{^C}D^{\alpha}_{0^{+}}f_{2}$ exists at $t=1$, and
$\,\leftidx{^C}D^{\alpha}_{0^{+}}f_{2}(1)= \frac{-1}{\Gamma(1-\alpha)}
\int_{0}^{1}\frac{d\xi}{(1-\xi)^{\alpha}}=
\frac{-1}{(1-\alpha)\Gamma(1-\alpha)}$.
\end{exm}
\end{rem}
\begin{rem}\label{remdfn9}
Pointwise Multiplication of the derivative $f^{\prime}$
of a function $f$ defined on $[a,\infty[$ by
the function $(t-a)^{1-\alpha}$ does not give the
physical properties we hope from
a fractional derivative. We show this by comparing $T^{a}_{\alpha}f(t)=(t-a)^{1-\alpha}f^{\prime}$ to
the Riemann-Liouville and Caputo fractional derivatives
for the sine and hyperbolic sine functions.
Similar differences show up with the cosine and hyperbolic cosine functions. Notice here that the Riemann-Liouville
fractional derivative coincides with the Caputo derivative for each of these functions.
We see the great difference in behaviour between $T^{a}_{\alpha}$ and the classical fractional
operators:
\begin{exm}
Let $g_{1}(t)=\sin{t}$. Then
$\left(T^{0}_{\alpha}g_{1}\right)(t)=t^{1-\alpha}
\cos{t}$. We can calculate
\begin{equation*}
\left(D^{\alpha}_{0^{+}}g_{1}\right)(t)=
\left(\leftidx{^C}D^{\alpha}_{0^{+}}g_{1}\right)(t) =
\frac{1}{\Gamma(1-\alpha)}\int_{0}^{t}
\frac{\cos{(t-\xi)}}{\xi^{\alpha}}d\xi.
\end{equation*}
Notice that $\left(T^{0}_{\alpha}g_{1}\right)(t)$ grows unboundedly with $t$. Contrarily, the fractional derivatives
$D^{\alpha}_{0^{+}}g_{1}$ and $\leftidx{^C}D^{\alpha}_{0^{+}}g_{1}$ are bounded. To see this, let $t>1$.
We have
\begin{equation}\label{bn1}
\left|\int_{0}^{1}
\frac{\cos{(t-\xi)}}{\xi^{\alpha}}d\xi\right|\leq
\int_{0}^{1}
\frac{1}{\xi^{\alpha}}d\xi=\frac{1}{1-\alpha}.
\end{equation}
Also, integrating by parts,
\begin{equation}\label{bn2}
\int_{1}^{t}
\frac{\cos{(t-\xi)}}{\xi^{\alpha}}d\xi=
\sin(t-1)-\alpha\int_{1}^{t}
\frac{\sin{(t-\xi)}}{\xi^{1+\alpha}}d\xi,
\end{equation}
and we have
\begin{equation}\label{bn3}
\left|\int_{1}^{t}
\frac{\sin{(t-\xi)}}{\xi^{1+\alpha}}d\xi\right|
\leq \int_{1}^{t}
\frac{1}{\xi^{1+\alpha}}d\xi=
\frac{1}{\alpha}\left(1-\frac{1}{t^{\alpha}}\right)<
\frac{1}{\alpha}.
\end{equation}
The boundedness of $D^{\alpha}_{0^{+}}g_{1}$, $\leftidx{^C}D^{\alpha}_{0^{+}}g_{1}$
follows from (\ref{bn1}), (\ref{bn2}), and (\ref{bn3}).
See Figure \ref{fig01}.
\begin{center}
\begin{figure}[H]
  \centering
  \includegraphics[width=10 cm]{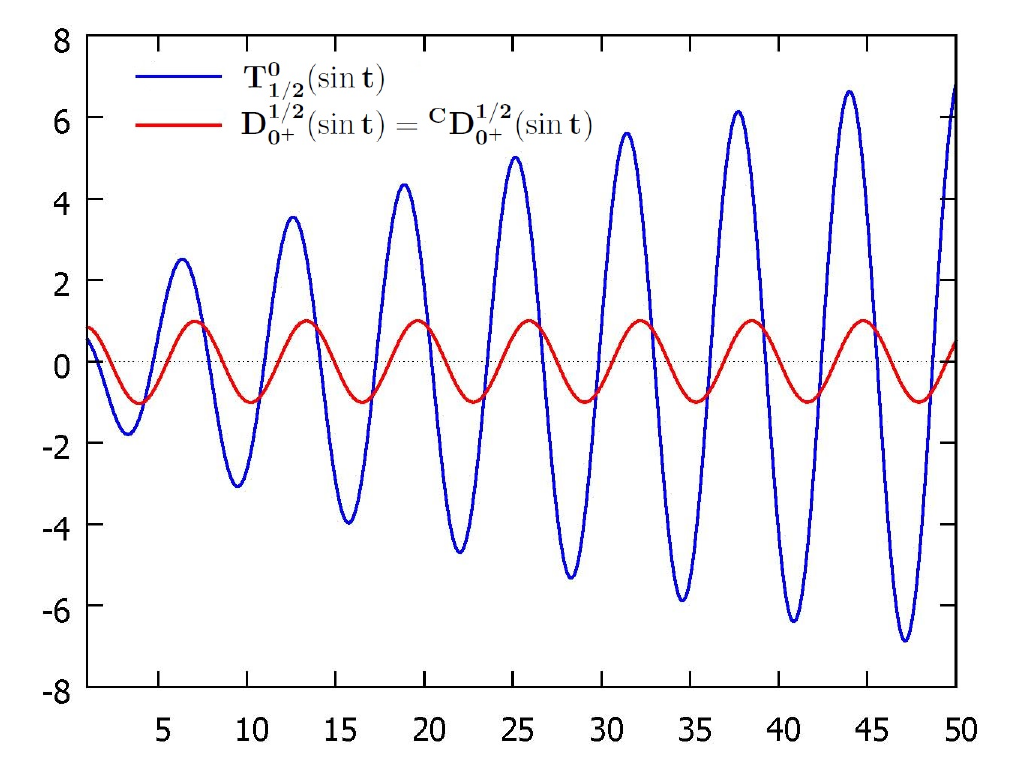}\\
  \caption{The behavior of  $T^{0}_{\alpha}g_{1}$ is very different from
  that of $D^{\alpha}_{0^{+}}g_{1}$, ${^C}D^{\alpha}_{0^{+}}g_{1}$}
  \label{fig01}
\end{figure}
\end{center}
\end{exm}
\begin{exm}
Let $g_{2}(t)=\sinh{t}$.
Then $\left(T^{0}_{\alpha}g_{2}\right)(t)=t^{1-\alpha}
\cosh{t}$.
On the other hand
\begin{equation*}
\left(D^{\alpha}_{0^{+}}g_{2}\right)(t)=
\left(\leftidx{^C}D^{\alpha}_{0^{+}}g_{1}\right)(t) =
\frac{1}{\Gamma(1-\alpha)}\int_{0}^{t}
\frac{\cosh{(t-\xi)}}{\xi^{\alpha}}d\xi.
\end{equation*}
The function $T^{0}_{\alpha}g_{2}$ grows much faster
than the fractional derivative. To prove this, we compute
\begin{align*}
&\lim_{t\rightarrow\infty}\frac{\int_{0}^{t}
\frac{\cosh{(t-\xi)}}{\xi^{\alpha}}d\xi}{t^{1-\alpha}
\cosh{t}}=
\lim_{t\rightarrow\infty}\frac{\int_{0}^{t}
\frac{1}{\xi^{\alpha}}{\left(
\cosh{\xi}-\tanh{t}\sinh{\xi}\right)}
d\xi}{t^{1-\alpha}}\\
=&\frac{1}{{1-\alpha}}\left(
\lim_{t\rightarrow\infty}
\frac{1}{\cosh{t}}-
\lim_{t\rightarrow\infty}\frac{\int_{0}^{t}
\frac{\sinh{\xi}}{\xi^{\alpha}}
d\xi}{t^{-\alpha}\cosh^{2}{t}}\right)=
\frac{-1}{{1-\alpha}}
\lim_{t\rightarrow\infty}\frac{\int_{0}^{t}
\frac{\sinh{\xi}}{\xi^{\alpha}}
d\xi}{t^{-\alpha}\cosh^{2}{t}}\\
=&\frac{-1}{{1-\alpha}}\lim_{t\rightarrow\infty}
\frac{1}{-\alpha\frac{\cosh^{2}{t}}
{t\sinh{t}}+2\cosh{t}}=0.
\end{align*}
See Figure \ref{fig02}.
\begin{center}
\begin{figure}[H]
  \centering
  \includegraphics[width=10 cm]{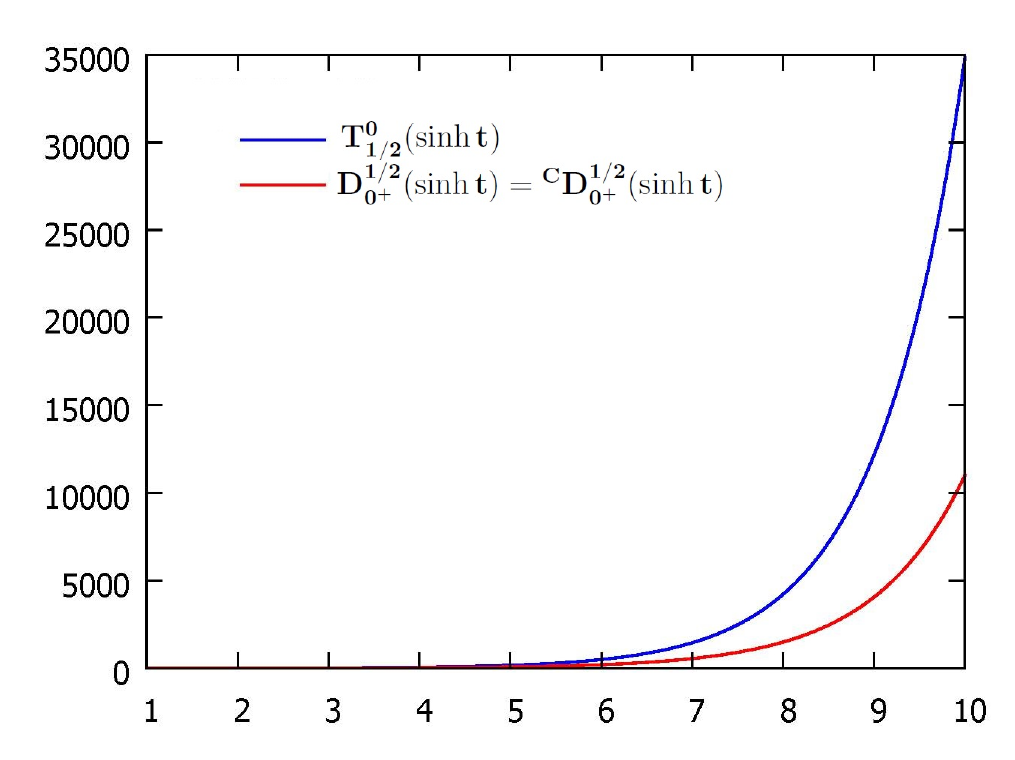}\\
  \caption{The behavior of  $T^{0}_{\alpha}g_{2}$ is very different from
  that of $D^{\alpha}_{0^{+}}g_{2}$, ${^C}D^{\alpha}_{0^{+}}g_{2}$}
  \label{fig02}
\end{figure}
\end{center}
\end{exm}
\end{rem}
\indent
Remarks \ref{remdfn} through \ref{remdfn9}
show the insignificance of the
results in \cite{Abdeljawad}. We illustrate how
the proofs presented in \cite{Abdeljawad}
reduce to trivial exercises of calculus. For example
\begin{thm}\label{thmj1}(\cite{Abdeljawad}, Theorem 2.11)
Assume $f,g:]a,\infty[\rightarrow \mathbb{R}$
are (left) $\alpha$-differentiable functions,
where $0<\alpha\leq 1$. Let $h(t)=f(g(t)$.
Then $h$ is (left) $\alpha$-differentiable, and
for all $t>a$ such that $g(t)\neq0$ we have
\begin{equation}\label{wrngcnc}
(T^{a}_{\alpha}h)(t)=
(T^{a}_{\alpha}f)(g(t)).
(T^{a}_{\alpha}g)(t).g(t)^{\alpha-1}.
\end{equation}
\end{thm}
\indent First of all, the conclusion (\ref{wrngcnc})
of Theorem \ref{thmj1} is incorrect.
It is correct if $a=0$. We consider this case.\\
\indent As noted in Remark \ref{remdfn}, if
$f$, $g$ are (left) $\alpha$-differentiable
on $]a,\infty[$,
then they are actually differentiable on $]a,\infty[$.
Moreover, by the identities (\ref{0qq3}),
\begin{align*}
(T^{0}_{\alpha}f)(g(t)).
(T^{0}_{\alpha}g)(t).g(t)^{\alpha-1}&=
g(t)^{1-\alpha}.f^{\prime}(g(t)).t^{1-\alpha}
g^{\prime}(t).g(t)^{\alpha-1}\\
&=
t^{1-\alpha}f^{\prime}(g(t)).g^{\prime}(t)\\
&=
t^{1-\alpha}(f(g(t)))^{\prime}=
(T^{0}_{\alpha}h)(t).
\end{align*}
\indent
Both the statement
and proof of the next theorem (\cite{Abdeljawad}, Theorem 4.1) we investigate
are incorrect. This implies that
Proposition 4.2 and Examples 4.1 through 4.3
in \cite{Abdeljawad} are also incorrect.
\begin{thm}\label{thmj2}(\cite{Abdeljawad}, Theorem 4.1)
Assume $f$ is an infinitely $\alpha$-differentiable function, for some $0<\alpha<1$ at a neighborhood of a point $t_{0}$. Then $f$
has the fractional power series expansion:
\begin{equation}\label{incrctxpn}
f(t)=\sum_{k=0}^{\infty}
\frac{(T^{t_{0}}_{\alpha}f)^{(k)}(t_{0})
(t-t_{0})^{k\alpha }}{{\alpha^{k}}{k!}},\;
t_{0}<t<t_{0}+R^{\frac{1}{\alpha}},\, R>0,
\end{equation}
where $(T^{t_{0}}_{\alpha}f)^{(k)}$
means the application of $T^{t_{0}}_{\alpha}$
$k$ times.
\end{thm}
\indent We will momentarily see that
the proof presented in \cite{Abdeljawad} is incorrect and the series (\ref{incrctxpn}) does not make sense.
The reason is $(T^{t_{0}}_{\alpha}f)^{(k)}(t_{0})$ is not well-defined unless $\alpha\leq 1/k$.\\
\indent The proof in \cite{Abdeljawad} begins with writing
\begin{equation}\label{prf1}
f(t)=c_{0}+c_{1}(t-t_{0})^{\alpha}+c_{2}(t-t_{0})^{2\alpha}
+c_{3}(t-t_{0})^{3\alpha}+...,
\end{equation}
and proceeds by applying $T^{t_{0}}_{\alpha}$
to both sides of (\ref{prf1}), then evaluating both sides at $t_{0}$, and repeating the process $k$ times. The coefficients
$c_{k}$ are inaccurately calculated:
\begin{equation*}
c_{k}=\frac{(T^{t_{0}}_{\alpha}f)^{(k)}(t_{0})
}{{\alpha^{k}}{k!}}.
\end{equation*}
\indent By Remark \ref{remdfn},
the assumption that $f$ is an infinitely $\alpha$-differentiable is equivalent to
assuming $f$ is infinitely differentiable.
Using (\ref{0qq3}), we get
\begin{eqnarray*}
   (T^{t_{0}}_{\alpha}f)^{(1)}(t)&=&
   (t-t_{0})^{1-\alpha}f^{\prime}(t),  \\
   (T^{t_{0}}_{\alpha}f)^{(2)}(t)&=&
  (1-\alpha)(t-t_{0})^{1-2\alpha}f^{\prime}(t)+ (t-t_{0})^{2-2\alpha}f^{\prime\prime}(t),\\
(T^{t_{0}}_{\alpha}f)^{(3)}(t)&=&
 (1-\alpha)(1-2\alpha)(t-t_{0})^{1-3\alpha}f^{\prime}(t)+ 3(1-\alpha)(t-t_{0})^{2-3\alpha}f^{\prime\prime}(t)+\\
&&+(t-t_{0})^{3-3\alpha}f^{\prime\prime\prime}(t),\\
(T^{t_{0}}_{\alpha}f)^{(4)}(t)&=&
 (1-\alpha)(1-2\alpha)(1-3\alpha)
 (t-t_{0})^{1-4\alpha}f^{\prime}(t)+ \\
 &&+(1-\alpha)(3(1-\alpha)+4(1-2\alpha))
(t-t_{0})^{2-4\alpha}f^{\prime\prime}(t)+\\
&&+6(1-\alpha)(t-t_{0})^{3-4\alpha}f^{\prime\prime\prime}(t)+
+(t-t_{0})^{4-4\alpha}f^{(4)}(t),\\
&.&\\
&.&\\
&.&\\
 (T^{t_{0}}_{\alpha}f)^{(k)}(t)&=&
 \prod_{j=1}^{k-1}
 (1-j\alpha)(t-t_{0})^{1-k\alpha}f^{\prime}(t)
+\sum_{j=2}^{k}a_{j,k}{(\alpha)}
 (t-t_{0})^{j-k\alpha}f^{(j)}(t),\, k\geq 2,
\end{eqnarray*}
where $a_{j,k}(\alpha)$ are constants that depend only
on $\alpha$.\\
\indent Since $\alpha \in\,]0,1[$, then there exists
$k_{0}\geq 2$ such that
$\frac{1}{k_{0}}<\alpha<\frac{2}{k_{0}}$.
Since, by the differentiability assumption,
the functions
$f^{(k)}$, $k\geq 0$, are all continuous,
then
\begin{equation*}
\lim_{t\rightarrow {t_{0}}^{+}}
\sum_{j=2}^{k_{0}}a_{j,k_{0}}{(\alpha)}
 (t-t_{0})^{j-k_{0}\alpha}f^{(j)}(t)=0.
\end{equation*}
However,
$(T^{t_{0}}_{\alpha}f)^{(k_{0})}(t)$ is ill-defined
at $t_{0}$ due to the term
$(t-t_{0})^{1-k_{0}\alpha}f^{\prime}(t)$.
So, unless $f$ is constant on a right neighborhood
of $t_{0}$ so that $f^{\prime}$ is identically zero
thereat, the coefficient $c_{k_{0}}$
in the series (\ref{prf1}) is at best indeterminate.
In fact, if the integer $j\geq 2$ is such that
$j>\frac{1}{\alpha}$, then none of the quantities
$(T^{t_{0}}_{\alpha}f)^{(k)}(t_{0})$, $k\geq j$,
is well-defined, because of the terms
$(t-t_{0})^{1-k\alpha}f^{\prime}(t)$.
This makes it impossible to compute any coefficient $c_{k}$ with $k\geq j$. Therefore, the series
(\ref{prf1}) does not make sense.\\
\indent  The examples, Examples 4.1 through 4.3, in \cite{Abdeljawad} that apply the incorrect expansion (\ref{prf1}) are conveniently for functions of the form
$f(t)=g\left(\left(\frac{t-t_{0}}{\alpha}\right)^{\alpha}\right)$. They seem to work because $(T^{a}_{\alpha}f)^{(k)}(t)=
g^{(k)}((t-t_{0})^{\alpha})$. The series (\ref{prf1}) fails even for the simplest analytic function $e^{t}$, or for the function $h(t)=e^{t^{\alpha}}$.
Let $t_{0}>a$, $t_{0}\neq 0$. Verifiably,  if $\alpha>\frac{1}{2}$, then $(T^{a}_{\alpha}h)^{(k)}(t_{0})$, $k\geq 2$, do not exit.
If $\alpha>\frac{1}{3}$, then
$(T^{a}_{\alpha}h)^{(k)}(t_{0})$, $k\geq 3$, do not exit. Generally, if $\alpha>\frac{1}{n}$ for some $n\geq 2$, then
$(T^{a}_{\alpha}h)^{(k)}(t_{0})$, $k\geq n$, do not exit. \\\\
\textbf{References}

\end{document}